\documentclass[12pt]{article}
\usepackage{amsmath,amssymb,amsthm}

\newtheorem{theorem}{Theorem}
\newtheorem{lemma}[theorem]{Lemma}

\newtheorem{corollary}[theorem]{Corollary}

\begin{document}
\title{Walks and paths in trees}
\author{B\'ela Bollob\'as\thanks{Department of Pure Mathematics and
Mathematical Statistics, University of Cambridge, UK}
\thanks{Research supported in part by NSF grants CCR-0225610,
DMS-0505550 and W911NF-06-1-0076} and
Mykhaylo Tyomkyn$^*$}
\maketitle

\begin{abstract}
Recently Csikv\'ari~\cite{csik} proved a conjecture of Nikiforov
concerning the number of closed walks on trees. Our aim is to extend
his theorem to all walks. In addition, we give a simpler
proof of Csikv\'ari's result and answer one of his questions in the negative.
Finally we consider an analogous question for paths rather than walks.
\end{abstract}

\section{Introduction}
In 2007, Pe\~{n}a, Rada and Gutman~\cite{pena} published the conjecture
that the maximum of the so-called Estrada index (see \cite{estr})
on trees of order $n$
is attained on a star, and its minimum is attained on a path.
This prompted Nikiforov~\cite{niki} to propose the stronger conjecture that
for a fixed value of $\ell$ the number of closed walks of length
$\ell$ on trees of order $n$ attains its extreme values on the
same graphs. Recently, Nikiforov's conjecture was proved by
Csikv\'ari~\cite{csik} by making use of a certain operation
inspired by a graph transformation which
Kelmans  (see \cite{kelm76a, kelm76b, kelm81}) defined in 1976
in order to prove
some results about the number of spanning trees of graphs.
In honour of Kelmans and Csikv\'ari, we call this tree-transformation the
{\em KC-transformation}. To define it,
let $x$ and $y$ be two vertices of a tree
$T$ such that every interior vertex of the unique $x$-$y$ path $P$ in $T$
has degree two, and write $z$ for the neighbour of $y$ on this path.
As usual, denote by $\Gamma(v)$ the set of neighbours of a vertex $v$.
The {\em KC-transform}  $KC(T,x,y)$
of the tree $T$ with respect to the path $P$ is
obtained from $T$ by deleting all edges
between $y$ and $\Gamma(y)\setminus z$ and adding the edges between $x$
and $\Gamma(y)\setminus z$ instead.
Note that $KC(T,x,y)$ and $KC(T,y,x)$ are isomorphic, so we may write
$KC(G,P)$ for this transform, without indicating in which `direction'
we take it.

Csikv\'ari proved that the
KC-transformation gives rise to a levelled poset of trees on $n$ vertices
with the star as the largest and the path as the smallest element.
To prove Nikiforov's conjecture, Csikv\'ari showed that, for any fixed
value of $\ell$, this transformation increases the number of {\em closed}
walks of length $\ell$.

In this paper, we extend this result to all walks of a given
length. Along the way, we give a considerably simpler proof of
Csikv\'ari's theorem, and answer a question he posed in the negative.

\vspace{5pt}
Let us remark that the analogous questions concerning paths rather
than walks in general graphs not only trees have been studied
since 1971; see, e.g., \cite{ahl-kat78,
alon81, alon86,  bb-erd, bb-sar1, bb-sar2, byer, fur92, katz71} and the references therein.

\vspace{5pt}
However, to the best of our knowledge, there has been no proper study of
the maximal number of paths of certain length in trees:
the final section of this paper is devoted to this topic.

\section{Further notation and Terminology}
Let $T$ be a tree and $T'=KC(T,p_0,p_k)$ its KC-transform along a path
$P=p_0p_1\dots p_k$, with edges $c_1,c_2, \dots , c_k$, where
$c_i=p_{i-1} p_i$.

Let $A$ and $B$ be the components of $p_0$ and $p_k$ in the graph
$T-E(P)$ obtained from $T$ by deleting the edges of $P$.
(Note that $A\cap B=\emptyset$ since $T$ is a tree.) Let $a_1, a_2, \dots$
be an enumeration of the edges of $A$, and $b_1, b_2, \dots$
an enumeration of the edges of $B$.
Label the edges of $T'$ in the same way, using the labels of deleted edges
between $p_k$ and $B$ for the corresponding new edges between $p_0$ and $B$.
As usual, write $G[X]$ for the edge-labelled subgraph of $G$ induced by a
subset $X$ of the vertices.

We encode the walks on $T$ and $T'$ by the sequences of traversed edges,
omitting the directions; in this way we assign to each walk a word on
the alphabet of the edge-labels $a_i,b_i$ and $c_i$.
Note that each walk of length at least two corresponds to a unique
word, while each single-letter word encodes two walks.

Let $\Omega(G)$ denote the set of all words corresponding to the walks on an
edge-labelled graph $G$. (Usually $G$ will be a subgraph of $T$ or $T'$.)
We refer to words in $\Omega(G)$ as {\em $G$-words}.
Furthermore, let us write $\Omega(x,G)$, $\Omega(G,y)$ and $\Omega(x,G,y)$
for the  set of all words encoding walks starting at $x\in G$,
finishing at $y\in G$ or doing both.
The set of  all {\em closed} $G$-words is
$\overline{\Omega}(G)=\bigcup_{x\in G}\Omega(x,G,x)$.

Denote the length of a word $W$ by $|W|$.
Let us write $\Omega_{\ell}(G)$ for the set of all $G$-words of length
$\ell$. Define $\Omega_{\ell}(x,G)$, $\Omega_{\ell}(G,y)$,
$\Omega_{\ell}(x,G,y)$ and $\overline{\Omega_{\ell}}(G)$ similarly.
Let $\omega(\hspace{1pt} . \hspace{2pt} )$
denote the size of the corresponding set
$\Omega(\hspace{2pt} . \hspace{2pt} )$,
for example, $\overline{\omega_{\ell}}(G)$ is the
size of $\overline{\Omega_{\ell}}(G)$.

\section{Block Structure}
In this section we fix a path $P=p_0 \dots p_k$ in our tree $T$, and write
$T'$ for $KC(T, P)$.
Every word on our alphabet can be decomposed into {\em blocks} of type $a$, $b$
and $c$, where an {\em $a$-block} is taken to be a maximal sequence of letters
of types $a$ and $c$ beginning and ending with an $a$. Likewise,
a {\em $b$-block} is a maximal
sequence of letters of types $b$ and $c$ beginning and ending with a
letter $b$. Finally,
and a {\em $c$-block} consists of all consecutive letters of type $c$
which are not in an $a$-block or a $b$-block.
We denote these blocks by the appropriate script letters.
It is always assumed that $\mathcal{A}$ and $\mathcal{B}$ are separated by $\mathcal{C}$, allowing $\mathcal{C}$
to be empty. Thus a typical block decomposition looks like
$\mathcal{C}\mathcal{A}\mathcal{C}\mathcal{B}\mathcal{C}\mathcal{A}\mathcal{C}\mathcal{B}\mathcal{C}\mathcal{A}$.
Let us call a block of a word {\em proper} if it is not the first or the last block of the word.

We would like to give a set of ``grammatical'' rules for blocks of each type in $T$-words and in $T'$-words.
In $W\in \Omega(T)$ every $a$-block $\mathcal{A}$ must be a $T[A\cup P]$-word. Moreover, we must have
$\mathcal{A}\in \Omega(p_0,T[A\cup P])$ if $\mathcal{A}$ is not the first block of $W$ and
$\mathcal{A}\in \Omega(T[A\cup P],p_0)$ if it is not the last one.
In particular, if $\mathcal{A}$ is proper, then $\mathcal{A} \in \Omega(p_0,T[A\cup P],p_0)$.
The same holds for $b$-blocks with $T[B\cup P]$ and $p_k$ in place of
$T[A\cup P]$ and $p_0$. A $c$-block $\mathcal{C}$ of $W$ must be a $P$-word, which satisfies
$\mathcal{C} \in \Omega(p_0,P)$ or $\mathcal{C} \in \Omega(p_k,P)$ if preceded by $\mathcal{A}$ or
$\mathcal{B}$ respectively, and $\mathcal{C} \in \Omega(P,p_0)$ or
$\mathcal{C} \in \Omega(P,p_k)$ if succeeded by $\mathcal{A}$ or $\mathcal{B}$ respectively.
In particular, if $\mathcal{C}$ is proper, then $\mathcal{C}\in \Omega(p_0,P,p_k)$
or $\mathcal{C}\in \Omega(p_k,P,p_0)$, depending on whether $\mathcal{C}$ is the middle block in the sequence
$\mathcal{A}\mathcal{C}\mathcal{B}$ or in $\mathcal{B}\mathcal{C}\mathcal{A}$. It is a trivial check that this set of rules is
complete, i.e.~if all blocks of a word $W$ satisfy them, then $W$ does indeed encode a walk on $T$.

Similar rules hold for a $T'$-word $W'$. More precisely, an $a$-block $\mathcal{A}$ has
to be a $T'[A\cup P]$-word, which lies in
$\Omega(p_0,T'[A\cup P])$ if $\mathcal{A}$ is not the first block of $W'$
and in $\Omega(T'[A\cup P],p_0)$ if it is not the last block of $W'$.
Likewise, a $b$-block $\mathcal{B}$ is a $T'[B\cup P]$-word with the properties
$\mathcal{B}\in \Omega(p_0,T'[B\cup P])$ if $\mathcal{B}$ is not at the beginning of $W'$
and $\mathcal{B} \in \Omega(T'[B\cup P], p_0)$ if it is not at the end.
Hence, proper blocks $\mathcal{A}$ and $\mathcal{B}$ must satisfy $\mathcal{A}\in \Omega(p_0,T'[A\cup P],p_0)$
and $\mathcal{B} \in \Omega(p_0,T'[B\cup P],p_0)$. Finally, a $c$-block of $W'$ must be a
$P$-word, satisfying $\mathcal{C} \in \Omega(p_0,P)$ if $\mathcal{C}$ is not the first block
of $W$ and $\mathcal{C} \in \Omega(P,p_0)$ if it is not the last one.
In particular, $\mathcal{C}\in \Omega(p_0,P,p_0)$ if $\mathcal{C}$ is proper.
As in the case of $T$-words, the above set of rules gives a complete characterization of blocks in $T'$-words.

Note that the edge-labelled graphs $T[A\cup P]$ and $T'[A\cup P]$ are identical. Hence, the sets of rules for an $a$-block in $T$ and $T'$ are the same.
To put it differently, any (first/last/proper) $a$-block $\mathcal{A}$ of $W\in \Omega(T)$ can be taken as a corresponding $a$-block in some $W'\in\Omega(T')$ and vice versa.
On the other hand, the edge-labelled graphs $T[B\cup P]$ and $T'[B\cup P]$ use the opposite numeration for the edges of $P$.
Therefore, the $b$-blocks in $T$-words and $T'$-words are the same up to replacing each $c_i$ with $c_{k+1-i}$. This operation will be called {\em conjugation} and, as usual, we denote the {\em conjugate} of $X$ by $\overline{X}$; trivially $\overline{\overline{X}}=X$. Thus, if $\mathcal{B}$ is a (first/last/proper) block of $W \in \Omega(T)$, then $\overline{\mathcal{B}}$ can be used as a corresponding $b$-block in some $W' \in \Omega(T')$ and vice versa.

\section{Closed walks}
In preparation for proving that the stars and the paths are extremal
for all walks, we give a new proof of Csikv\'ari's theorem.
\begin{theorem} \label{thm:csik}
For every $\ell \ge 1$, the KC-transformation on trees increases
the number of closed walks of length $\ell$.
\end{theorem}
Needless to say, `increases' is used in the usual weak sense that
the number does not become strictly smaller.
As in \cite{csik}, this theorem has the following consequence.
\begin{corollary}\label{cor:csik}
The number of closed walks of length $\ell$ in a tree $T$ on $n$ vertices
is maximal when $T$ is a star and minimal when $T$ is a path.
\end{corollary}

\vspace{5pt}
We shall prove Theorem~\ref{thm:csik} by
defining an injective mapping
$f\colon \overline{\Omega_{\ell}}(T)\rightarrow \overline{\Omega_{\ell}}(T')$.
Given a word $W \in \overline{\Omega_{\ell}}(T)$, our definition of $f(W)$
will depend on the `type' of $W$. Indeed, it is easily seen that
$W$ is precisely one of
the following (mutually exclusive) five types.

\vspace{5pt}
\textbf{0.}
$W$ is a single $c$-block.

\vspace{5pt}
\textbf{1.1.}
$W$ has an even number number of proper $\mathcal{C}$'s and the first $a$ appears before the first $b$. In other words, $W$ can be written as \[W = (\mathcal{C}_0)\mathcal{A}_1\mathcal{C}_1\mathcal{B}_2\mathcal{C}_2\mathcal{A}_3\dots \mathcal{C}_{2m-1}\mathcal{B}_{2m}\mathcal{C}_{2m}\mathcal{A}_{2m+1}(\mathcal{C}_{2m+1}).\]

\vspace{5pt}
\textbf{1.2.}
$W$ has an even number of proper $\mathcal{C}$'s and the first $b$ appears before the first $a$, i.e. \[W=(\mathcal{C}_0)\mathcal{B}_1\mathcal{C}_1\mathcal{A}_2\mathcal{C}_2\mathcal{B}_3\dots \mathcal{C}_{2m-1}\mathcal{A}_{2m}\mathcal{C}_{2m}\mathcal{B}_{2m+1}(\mathcal{C}_{2m+1}).\]

\vspace{5pt}
\textbf{2.1.}
$W$ has an odd number of proper $\mathcal{C}$'s and the first $a$ appears before the first $b$, i.e. \[W=(\mathcal{C}_0)\mathcal{A}_1\mathcal{C}_1\mathcal{B}_2\mathcal{C}_2\dots \mathcal{A}_{2m-1}\mathcal{C}_{2m-1}\mathcal{B}_{2m}(\mathcal{C}_{2m}).\]

\vspace{5pt}
\textbf{2.2.}
$W$ has an odd number of proper $\mathcal{C}$'s and the first $b$ appears before the first $a$, i.e. \[W=(\mathcal{C}_0)\mathcal{B}_1\mathcal{C}_1\mathcal{A}_2\mathcal{C}_2\dots \mathcal{B}_{2m-1}\mathcal{C}_{2m-1}\mathcal{A}_{2m}(\mathcal{C}_{2m}).\]

\vspace{5pt}
The first and the last $c$-blocks are taken into parentheses, in order to
indicate that they might not exist. We should like to define our map
$f$ such that it maps $W$ to $W'\in \overline{\Omega_{\ell}}(T')$ of the same type.
In that way we ensure that the images of the types are disjoint, and therefore
$f$ is injective on the whole of its domain.
Furthermore, in the first three cases we actually construct a more general
injective mapping $f\colon \Omega_{\ell}(T)\rightarrow \Omega_{\ell}(T')$ which
happens to map closed words to closed $T'$-words.
This fact will be useful in the next section.

\vspace{5pt}
\textbf{Case 0.}
Set $f(W)=W$.

\vspace{5pt}
\textbf{Case 1.1.}
Take $\mathcal{C}_1\in \Omega(p_0,P,p_k)$ and split it at its walk's last point of visit to $p_0$ into $\mathcal{C}_{1,1} \in \Omega(p_0,P,p_0)$, possibly empty, and $\mathcal{C}_{1,2}\in \Omega(p_0,P,p_k)$. Doing the same to each proper $\mathcal{C}_{2i+1}$, we can write  \[W=(\mathcal{C}_0)\mathcal{A}_1\mathcal{C}_{1,1}\mathcal{C}_{1,2}\mathcal{B}_2\mathcal{C}_2\mathcal{A}_3\dots \mathcal{C}_{2m-1,1}\mathcal{C}_{2m-1,2}\mathcal{B}_{2m}\mathcal{C}_{2m}\mathcal{A}_{2m+1}(\mathcal{C}_{2m+1}).\] Now define $f(W)$ to be the word \[W'=(\mathcal{C}_0)\mathcal{A}_1\mathcal{C}_{1,1}\overline{\mathcal{B}_2}\overline{\mathcal{C}_{1,2}^r}\mathcal{C}_2\mathcal{A}_3\dots \mathcal{C}_{2m-1,1}\overline{\mathcal{B}_{2m}}\overline{\mathcal{C}_{2m-1,2}^r}\mathcal{C}_{2m}\mathcal{A}_{2m+1}(\mathcal{C}_{2m+1}),\] where $X^r$ stands for $X$ spelt backwards (note that $(X^r)^r=X$ and $\overline{X^r}=\overline X^r$).

Since both the conjugation and the reversed spelling are length preserving, it follows that $|W'|=|W|=\ell$. One can easily convince oneself that the blocks of $W'$ are indeed as the above representation suggests, namely various $\mathcal{A}_i$, $\overline{\mathcal{B}_i}$, $\mathcal{C}_{2i-1,1}$, $\overline{\mathcal{C}_{2i-1,2}^r}\mathcal{C}_{2i}$ and perhaps $\mathcal{C}_0$ and $\mathcal{C}_{2m+1}$.~Thus, in order to show that $W'\in \Omega(T')$, we must check that these blocks meet the corresponding conditions.

For $\mathcal{A}_i$ and $\overline{\mathcal{B}_i}$, $\mathcal{C}_0$ and $\mathcal{C}_{2m+1}$ this follows from the observations in the previous section.
The fact $\mathcal{C}_{2i-1,1}\in \Omega(p_0,P,p_0)$ follows from the definition of $\mathcal{C}_{2i-1,1}$. It follows also from definitions that $\mathcal{C}_{2i-1,2}\in \Omega(p_0,P,p_k)$. This implies $\mathcal{C}_{2i-1,2}^r\in \Omega(p_k,P,p_0)$ and hence $\overline{\mathcal{C}_{2i-1,2}^r}\in \Omega(p_0,P,p_k)$. Therefore, since $\mathcal{C}_{2i}\in \Omega(p_k,P,p_0)$, we have $\overline{\mathcal{C}_{2i-1,2}^r}\mathcal{C}_{2i}\in \Omega(p_0,P,p_0)$. This proves that $W'\in \Omega_{\ell}(T')$.

Observe that $\overline{\mathcal{C}_{2i-1,2}^r}$ and $\mathcal{C}_{2i}$ can be recovered from  $\overline{\mathcal{C}_{2i-1,2}^r}\mathcal{C}_{2i}$ by splitting the latter at the first point of visit of its walk to $p_k$, that is after the first occurrence of the letter $c_k$. This fact implies that the mapping $f$ is injective, as we can define an inverse mapping by sending \[W'=(\mathcal{C}_0)\mathcal{A}_1\mathcal{C}_1\mathcal{B}_2\mathcal{C}_2\mathcal{A}_3\dots \mathcal{C}_{2m-1}\mathcal{B}_{2m}\mathcal{C}_{2m}\mathcal{A}_{2m+1}(\mathcal{C}_{2m+1}) \in \Omega_\ell(T')\] to the word \[W=(\mathcal{C}_0)\mathcal{A}_1\mathcal{C}_1\overline{\mathcal{C}_{2,1}^r\mathcal{B}_2}\mathcal{C}_{2,2}\mathcal{A}_3\dots \mathcal{C}_{2m-1}\overline{\mathcal{C}_{2m,1}^r\mathcal{B}_{2m}}\mathcal{C}_{2m,2}\mathcal{A}_{2m+1}(\mathcal{C}_{2m+1}),\] where $\mathcal{C}_{2i,1}$ is the initial segment of $\mathcal{C}_{2i}$ up to the first visit to $p_k$, as was just mentioned.

Note that if $W$ is closed, then so is $f(W)$, since the walk encoded by the latter starts and ends at the same point as the walk of $W$.

\vspace{5pt}
\textbf{Case 1.2.}
Similarly to the previous case, let us split each proper $\mathcal{C}_{2i+1} \in \Omega(p_k,P,p_0)$ at its walk's last point of visit to $p_k$ into $\mathcal{C}_{2i+1,1}\in \Omega(p_k,P,p_k)$ and $\mathcal{C}_{2i+1,2}\in \Omega(p_k,P,p_0)$.
We obtain\[W=(\mathcal{C}_0)\mathcal{B}_1\mathcal{C}_{1,1}\mathcal{C}_{1,2}\mathcal{A}_2\mathcal{C}_2\mathcal{B}_3\dots \mathcal{C}_{2m-1,1}\mathcal{C}_{2m-1,2}\mathcal{A}_{2m}\mathcal{C}_{2m}\mathcal{B}_{2m+1}(\mathcal{C}_{2m+1}).\] Define $f(W)$ to be the word \[W'=\overline{(\mathcal{C}_0)\mathcal{B}_1\mathcal{C}_{1,1}}\mathcal{A}_2\mathcal{C}_{1,2}^{r}\overline{\mathcal{C}_2\mathcal{B}_3\mathcal{C}_{3,1}}\mathcal{A}_4\mathcal{C}_{3,2}^r\dots \mathcal{A}_{2m}\mathcal{C}_{2m-1,2}^{r}\overline{\mathcal{C}_{2m}\mathcal{B}_{2m+1}(\mathcal{C}_{2m+1})}.\]

Again, we have $|W'|={\ell}$, by the length-invariance of conjugation and reversion. By the observation in the previous section the blocks $\mathcal{A}_i$ and $\overline{\mathcal{B}_i}$ satisfy the rules for $a$-blocks and $b$-blocks in a $T'$-word.~Similarly as in case 1.1.~we verify that
$\overline{\mathcal{C}_{2i-1,1}} \in \Omega(p_0,P,p_0)$ since $\mathcal{C}_{2i-1,1}\in \Omega(p_k,P,p_k)$, and $\mathcal{C}_{2i-1,2}^r\overline{\mathcal{C}_{2i}} \in \Omega(p_0,P,p_0)$ since $\mathcal{C}_{2i-1,2}^r \in \Omega(p_0,P,p_k)$ and $\overline{\mathcal{C}_{2i}}\in \Omega(p_k,P,p_0)$. It should be also remarked that $\overline{\mathcal{C}_0}\in \Omega(P,p_0)$ and $\overline{\mathcal{C}_{2m+1}}\in \Omega(p_0,P)$ meet the conditions for a block in a $T'$-word as well. Therefore $W \in \Omega_{\ell}(T')$.

The inverse mapping can be defined as the function that takes \[W' = (\mathcal{C}_0)\mathcal{B}_1\mathcal{C}_1\mathcal{A}_2\mathcal{C}_2\mathcal{B}_3\dots \mathcal{C}_{2m-1}\mathcal{A}_{2m}\mathcal{C}_{2m}\mathcal{B}_{2m+1}(\mathcal{C}_{2m+1})\in \Omega_\ell(T')\] to the word \[W = \overline{(\mathcal{C}_0)\mathcal{B}_1\mathcal{C}_{1}}\mathcal{C}_{2,1}^{r}\mathcal{A}_2\overline{\mathcal{C}_{2,2}\mathcal{B}_3\mathcal{C}_3}\mathcal{C}_{4,1}^{r}\dots \mathcal{C}_{2m,1}^{r}\mathcal{A}_{2m}\overline{\mathcal{C}_{2m,2}\mathcal{B}_{2m+1}(\mathcal{C}_{2m+1})},\] where $\mathcal{C}_{2i,1}$ is the initial part of $\mathcal{C}_{2i}$ up to its walk's first visit to $p_k$. Hence, $f$ is injective.

Since the first and the last blocks of $W$ are conjugated to the first and last block of $W'$ respectively, the fact $W\in \overline{\Omega}(T)$ would imply $W'\in \overline{\Omega}(T')$.

\vspace{5pt}
Recall that each single-letter word $W$ encodes, as was remarked in section 2, precisely two different words. Since every such $W$ must have one of the types 0, 1.1. or 1.2., we can deduce from the definition of $f$ on these types that $f(W)=W$. Therefore, $f$ does indeed give rise to an injective mapping of walks on $T$ into walks on $T'$ of the same length.

\vspace{5pt}
\textbf{Case 2.1.}
Unlike in the previous two cases, where we managed to construct injective length preserving mappings from $T$-words into $T'$-words, which also happened to map $\overline{\Omega}(T)$ into $\overline{\Omega}(T')$, in this and the next case we confine ourselves to closed words. Here and in the next case we must have $W\in \Omega(p_i,T,p_i)$ for some $p_i \in P$, otherwise there would have been an even number of proper $\mathcal{C}$'s. Therefore,  $\mathcal{B}_{2m}\in \Omega(T,p_k)$ even if $\mathcal{C}_{2m}$ does not exist, i.e. if $\mathcal{B}_{2m}$ is not proper. Split each proper $\mathcal{C}_{2i+1}$ into $\mathcal{C}_{2i+1,1}$ and $\mathcal{C}_{2i+1,2}$ as in case 1.1:\[W=(\mathcal{C}_0)\mathcal{A}_1\mathcal{C}_{1,1}\mathcal{C}_{1,2}\mathcal{B}_2\mathcal{C}_2\dots \mathcal{A}_{2m-1}\mathcal{C}_{2m-1,1}\mathcal{C}_{2m-1,2}\mathcal{B}_{2m}(\mathcal{C}_{2m}).\] Define $f(W)$ to be \[W'=(\mathcal{C}_0)\mathcal{A}_1\mathcal{C}_{1,1}\overline{\mathcal{B}_2\mathcal{C}_{1,2}^r}\mathcal{C}_2\dots \mathcal{A}_{2m-1}\mathcal{C}_{2m-1,1}\overline{\mathcal{B}_{2m}\mathcal{C}_{2m-1,2}^r}(\mathcal{C}_{2m}).\]

As in the previous cases, we have $|W'|=\ell$. In order to see that $W \in \Omega(T')$, the only thing that needs to be checked beyond case 1.1.~is that $\overline{\mathcal{B}_{2m}} \in \Omega(P,p_0)$, since this block is proper in $W'$. This, however, follows from the above observation that $\mathcal{B}_{2m}\in \Omega(T,p_k)$. Hence, $W' \in \Omega(T')$. The injectivity of $f$ is provided by the same inverse mapping as in case 1.1.

To see that $W' \in \overline{\Omega}(T')$ note that the starting points of $W$ and $W'$ coincide, hence it is enough to show that the endpoints do as well. This is obvious when $\mathcal{C}_{2m}$ exists. If it does not, then by the virtue of $B_{2m}\in \Omega(T,p_k)$ and $\overline{\mathcal{C}_{2m-1,2}^r}\in \Omega(T,p_k)$ the walks of $W$ and $W'$ both end at $p_k$. It follows that $W' \in \overline{\Omega_{\ell}}(T')$.

\vspace{5pt}
\textbf{Case 2.2.}
Recall that, as in case 2.1., the walk of $W$ starts and finishes on $P$. So, whether $C_0$ and $C_{2m}$ exist or not, we have $\mathcal{B}_1 \in \Omega(p_k,T)$ and $\mathcal{A}_{2m}\in \Omega(T,p_0)$. Split each $\mathcal{C}_{2i+1}$ into $\mathcal{C}_{2i+1,1}$ and $\mathcal{C}_{2i+1,2}$ as in case 1.2: \[W=(\mathcal{C}_0)\mathcal{B}_1\mathcal{C}_{1,1}\mathcal{C}_{1,2}\mathcal{A}_2\mathcal{C}_2\dots \mathcal{B}_{2m-1}\mathcal{C}_{2m-1,1}\mathcal{C}_{2m-1,2}\mathcal{A}_{2m}(\mathcal{C}_{2m}).\] Define $f(W)$ as \[W'=\overline{(\mathcal{C}_0)\mathcal{B}_1\mathcal{C}_{1,1}}\mathcal{A}_2\mathcal{C}_{1,2}^{r}\overline{\mathcal{C}_2\mathcal{B}_3\mathcal{C}_{3,1}}\mathcal{A}_4\mathcal{C}_{3,2}^r\dots \mathcal{A}_{2m}\mathcal{C}_{2m-1,2}^{r}(\overline{\mathcal{C}_{2m}}).\]
$W'\in \Omega(T')$ follows from the observations from case 1.2 and the fact $\mathcal{A}_{2m}\in \Omega(T,p_0)$, proving that $\mathcal{A}_{2m}$ is a suitable proper $a$-block of a $T'$-word. The injectivity of $f$ follows as in case 1.2.

To show that $W'$ is closed, observe that since $\mathcal{B}_1 \in \Omega(p_k,T)$ and the first blocks of $W$ and $W'$ are conjugated, no matter if $\mathcal{C}_0$ exists or not, we have $W\in \Omega(p_i,T)$ and $W'\in \Omega(p_{k-i},T')$ for some $i$. On the other hand, since $\mathcal{A}_{2m}\in \Omega(T,p_0)$ and $\mathcal{C}_{2m-1,2}^r \in \Omega(T,p_k)$, we also have $W\in \Omega(T,p_i)$ and $W'\in \Omega(T,p_{k-i})$. Therefore, $W'$ is closed as well.

\section{General walks}
In this section we shall prove our main results about walks.

\begin{theorem}\label{thm: Main}
For every $\ell \ge 1$, the KC-transformation on trees increases
the number of closed walks of length $\ell$.
\end{theorem}
\begin{corollary}
The number of closed walks of length $\ell$ in a tree $T$ on $n$ vertices
is maximal when $T$ is a star and minimal when $T$ is a path.
\end{corollary}

Let us first assume that $k$ is even. We can prove the following useful fact about walks on $T[B\cup P]$.

\begin{lemma}\label{lemma2}
If $P$ is of even length, then for every $\ell\geq 1$ holds
\[\omega_{\ell}(p_0,T[B\cup P])-\omega_{\ell}(p_0,P)\leq \omega_{\ell}(p_k,T[B\cup P])-\omega_{\ell}(p_k,P).\]
\end{lemma}

In other words, among all walks on $T[B\cup P]$ that visit $B$, there are fewer walks starting at $p_0$ than at $p_k$. We would like to give an injective mapping $g$ from the former into the latter.

Take a word $W\in \Omega_{\ell}(p_0,T[B\cup P])\setminus \Omega(P)$. Since $W$ contains at least one $b$, its walk visits $p_{k/2}$. So we can decompose $W$ at the first point of visit to $p_{k/2}$ into $W_1\in \Omega(p_0,P,p_{k/2})$ and $W_2\in \Omega(p_{k/2},T[B\cup P])\setminus \Omega(P)$. Define $g(W)$ to be $W'=\overline{W_1}W_2$. It is immediate that $|W'|=|W|=\ell$. Since $\overline{W_1}\in \Omega(p_k, T[B\cup P],p_{k/2})$, the word $W'$ is in $\Omega_{\ell}(p_k, T[B\cup P])$. Since $W_2\notin \Omega(P)$, we have $W'\notin \Omega(P)$. There mapping $g$ is self-inverse, thus injective.

By identifying $\Omega_{\ell}(p_k,T[B\cup P])$ with $\Omega_{\ell}(p_0,T'[B\cup P])$ via conjugation we can rewrite the above statement as follows.

\begin{corollary}\label{corlemma2}
 If $P$ is of even length, then for every $\ell\geq 1$ holds
\[\omega_{\ell}(p_0,T[B\cup P])-\omega_{\ell}(p_0,P)\leq \omega_{\ell}(p_0,T'[B\cup P])-\omega_{\ell}(p_0,P).\]
\end{corollary}

\vspace{5pt}

For an odd $k$ we can prove an analogous inequality with $\ell-1$ instead of $\ell$ on the right hand side.

\begin{lemma}\label{lemma3}
If $P$ is of odd length and $B$ is not empty, then for any $\ell$ holds \[\omega_{\ell}(p_0,T[B\cup P])-\omega_{\ell}(p_0,P)\leq \omega_{\ell-1}(p_k,T[B\cup P])-\omega_{\ell-1}(p_k,P).\]
\end{lemma}

Since $\Omega_{\ell}(p_0,T[B\cup P])$ can be identified with $\Omega_{\ell-1}(p_1,T[B\cup P])$, it is enough to show that the latter contains fewer elements than $\Omega_{\ell-1}(p_k,T[B\cup P])\setminus \Omega(P)$. Let $u$ be a designated neighbour of $p_k$ in $B$ and let $P'$ be $T[P\cup\left\{u\right\}]$, i.e.~the extension of $P$ to $u$.~Note that $P'$ is a path of even length and $p_{(k+1)/2}$ is its midpoint. As in the proof of lemma \ref{lemma2} we can split $W\in \Omega_{\ell-1}(p_k,T[B\cup P])\setminus \Omega(P)$ at the first visit of its walk to $p_{(k+1)/2}$ into $W_1$ and $W_2$. Define $g(W)$ to be $W'=W'_1W_2$, where $W'$ is the image of $W_1$ under reflexion on $P'$. It follows as in the proof of lemma \ref{lemma2} that $W'\in \Omega_{\ell-1}(p_k,T[B\cup P])\setminus \Omega(P)$. Again, $g$ is self-inverse, therefore injective.

As with lemma \ref{lemma2}, we obtain the an immediate consequence for $T'[B\cup P]$. Since we can extend every $\Omega_{\ell-1}(p_k,T[B\cup P])\setminus \Omega(P)$-word by an appropriate letter to a word of length $\ell$, we can replace $\ell-1$ by $\ell$ on the right hand side and drop the requirement of $B$ being not empty.

\begin{corollary}\label{corlemma3}
If $P$ is of odd length, then for any $\ell$ holds \[\omega_{\ell}(p_0,T[B\cup P])-\omega_{\ell}(p_0,P)\leq \omega_{\ell}(p_0,T'[B\cup P])-\omega_{\ell}(p_0,P).\]
\end{corollary}

Let us summarize.

\begin{corollary}\label{cortotal}
For every $\ell\geq 1$ holds \[\omega_{\ell}(p_0,T[B\cup P])-\omega_{\ell}(p_0,P)\leq \omega_{\ell}(p_0,T'[B\cup P])-\omega_{\ell}(p_0,P).\]
\end{corollary}

\vspace{5pt}

It goes without saying that analogous statements hold in $T[A\cup P]$ and $T'[A\cup P]$.

\vspace{10pt}

\textbf{Proof of Theorem \ref{thm: Main}.} As in the case of closed walks, we would like to define an injective mapping $h\colon \Omega_{\ell}(T)\rightarrow \Omega_{\ell}(T')$. For this sake we use the classification of words according to their type, as defined in the previous section. Like for the closed walks, we would like $h$ to map words within their type. This fact and the injectivity of $h$ on each single type would ensure that $h$ is injective on its whole domain.

If $W \in \Omega_{\ell}(T)$ is of type 0, 1.1. or 1.2., i.e.~if $W$ has an even number of proper $c$-blocks, then define $h$ to be the mapping $f$ from previous section. Recall that if $W$ has one of the above types, then $f$ is well-defined for general walks, injective, length-preserving and maps $W$ to $W'\in \Omega_{\ell}(T')$ of the same type.

Suppose now that $W$ is of type 2.1. Let us write it in the generic block from \[W=(\mathcal{C}_0)\mathcal{A}_1\mathcal{C}_1\mathcal{B}_2\mathcal{C}_2\dots \mathcal{A}_{2m-1}\mathcal{C}_{2m-1}\mathcal{B}_{2m}(\mathcal{C}_{2m}).\] Decompose $W$ into $W_1=(\mathcal{C}_0)\mathcal{A}_1\mathcal{C}_1\mathcal{B}_2\mathcal{C}_2\dots \mathcal{A}_{2m-1}\in \Omega(T,p_0)$ and $W_2=\mathcal{C}_{2m-1}\mathcal{B}_{2m}(\mathcal{C}_{2m})\in \Omega(p_0,T[B\cup P])$. Note that $W_1$ has an even number of proper $\mathcal{C}$'s.

Define $h(W)$ to be $W' = f(W_1)g(W_2)$ where $f$ is as above and $g$ is an injective mapping provided by corollary \ref{cortotal}. This is a $T'$-word, since $f(W_1)\in \Omega(T',p_0)$, by construction of $f$ on type 1.1. and $g(W_2)\in \Omega(p_0,T')$ by definition of $g$.

Note that $W'$ is again of type 2.1, as  $g(W)\in \Omega(T'[B\cup P])\setminus \Omega(P)$. Since $f(W_1)$ and $g(W_2)$ can be recovered from $W'$ by splitting it in the same way as we split $W$ and both $f$ and $g$ are injective, the mapping $h$ is injective on this type as well. The length of $W'$ equals the length of $W$, since $f$ and $g$ are length-preserving.

Finally, suppose that $W$ is of type 2.2. \[W=(\mathcal{C}_0)\mathcal{B}_1\mathcal{C}_1\mathcal{A}_2\mathcal{C}_2\dots \mathcal{B}_{2m-1}\mathcal{C}_{2m-1}\mathcal{A}_{2m}(\mathcal{C}_{2m}).\] Decompose it as above into \[W_1=(\mathcal{C}_0)\mathcal{B}_1\mathcal{C}_1\mathcal{A}_2\mathcal{C}_2\dots \mathcal{B}_{2m-1}\] and \[W_2=\mathcal{C}_{2m-1}\mathcal{A}_{2m}(\mathcal{C}_{2m}).\] Applying the analogue of corollary \ref{cortotal} to $W_2$ and to $A\cup P$ instead of $B\cup P$, we can define $h(W)$ to be $W'=f(W_1)g(W_2)$ where $f$ is as in case 1.2. Similarly to the definition of $h$ in the previous case, $W'$ is a $T'$-word of type 2.2., and the mapping is length-preserving and injective.

Therefore, the total mapping $h$ is injective and length-preserving. This proves the theorem.

\section{An answer  to Csikv\'ari's question}
Csikv\'ari~\cite{csik} proved also that $D(T) = \sum_{x,y\in T}d(x,y)$ is
decreased and $W_{\ell}$, the number of closed walks of length $\ell$, is
increased by a KC-transformation. Based on this fact he asked whether
$D(T_1)<D(T_2)$ implies $W_{\ell}(T_1)>W_{\ell}(T_2)$ for any two trees $T_1$ and
$T_2$ of size $n$. Our next aim is to show that this is not the case.

Let $T_1$ be a ``broom'', i.e. a path on length $(2-c)k$ with $ck$ leaves attached to one of its endpoints, where $k$ is large and $c$ will be specified later. Let $T_2$ be a ``double broom'', i.e. a path of length $k$ with $k/2$ additional leaves attached to each endvertex. Both trees have about $2k$ vertices and edges --- since our estimates will be asymptotic, there is no need for exact counting. By chosing an appropriate $c$ we would like to achieve $D(T_1)>D(T_2)$ and $W_4(T_1)>W_4(T_2)$ simultaneously. Note that $\ell=4$ is the smallest non-trivial case.

The actual reason why this construction produces a counterexample is that the single-brooms on a given number of vertices form a totally ordered subset of the poset induced by $KC$-transformations. Indeed, each time we apply the $KC$-transformation to the centre of the star and the adjacent non-leaf, we obtain a new broom with the length of the path decreased by $1$ and the order of the star increased by $1$. Along these transformations $D(T_1)$ decreases and $W_\ell(T_1)$ increases in rather small steps from one extremal case (path) to the other (star). Therefore, a counterexample to  Csikv\'ari's question exists, unless for every tree $T$ on $n$ vertices the values of $D(T)$ and $W_\ell(T)$ lie between the respective values of two 'consecutive' brooms, which strongly suggests that $W_\ell$ is a function of $D$. The latter appears rather unlikely.

To construct an explicit counterexample, let us count ordered pairs of adjacent edges in $T_1$ and $T_2$. It is the same, up to a constant factor and a negligible error term, as counting closed walks of length $4$. Since $k$ is large, we can ignore the path and count just edge pairs in the stars. We obtain $c^2k^2+O(k)$ such pairs in $T_1$ and $k^2/2+O(k)$ in $T_2$.

Now let us estimate $D(T_1)$ and $D(T_2)$. For a path $P$ of length $t$ we get \[D(P)= \sum_{i=1}^{t}\frac{i(i+1)}{2} = \frac{1}{6}t^3+O(t^2).\] Therefore, in $T_1$ we count about $(2-c)^3k^3/6$ for the distances between two points on the path, about $ck[(2-c)k]^2/2$ for the distances between a point in the star and a point on the path and about $2(ck)^2/2$ for the distances inside the star. The last term is negligibly small, so we can write \[D(T_1)=\left(\frac{(2-c)^3}{6}+\frac{c(2-c)^2}{2} + o(1)\right)k^3 = \left(\frac{(2-c)^2(1+c)}{3}+ o(1)\right)k^3.\]

In a similar fashion we can estimate $D(T_2)$: we count about $k^3/6$ for the distances on the path, about $2(k/2)(k^2/2)$ for distances between a point in one of the stars and a point on the path and finally about $(k/2)^2k$ for the distances between two points in different stars. The distances inside each star contribute only about $2(k/2)^2$ and can thus be neglected. In total we obtain \[D(T_2)=\left(\frac{1}{6}+\frac{1}{2}+\frac{1}{4}+o(1)\right)k^3 = \left(\frac{11}{12}+o(1)\right)k^3.\]

Now we want to choose a $c$ such that $D(T_1)>D(T_2)$ and $W_4(T_1)>W_4(T_2)$, giving the negative answer to Csikv\'ari's question. By the above estimates we want $c$ to satisfy $c^2>1/2$ and $f(c)=(2-c)^2(1+c)>11/4$. The first condition can be stated as $c>\sqrt{1/2}\approx 0.707$. The second is met for $c<0.744$, hence any $c$ between $\sqrt{1/2}$ and $0.744$ would satisfy both conditions.

To be more concrete, one could take $c=0.72$ and $k=1000$. A straightforward calculation confirms that this is ineed a counterexample.

It can be shown in a very similar way that for each $\ell\geq 2$ there is a $c$ such that the above construction yields a counterexample for closed walks of length $2\ell$ or general walks of length $\ell$. One can easily convince oneself that for a large enough $k$ the only walks in $T_1$ and $T_2$ that make a significant contribution to the total number lie entirely in one of the stars. Hence, in order to construct a counterexample, we need a $c$ which satisfies $c^{\ell}>2(1/2)^\ell$ and $f(c)>11/4$. Notice that any $c$ that satisfies the first inequality for $\ell=2$ does so for all values of $\ell$, whereas the second inequality does not depend on $\ell$. Therefore the above $c=0.72$ satisfies both inequalities, and thus can be used to construct a counterexample for any arbitrary $\ell$ (however we would have to choose a different $k$ each time).

It is an interesting question, whether there is a ``universal'' counterexample for all values of $\ell$.

\section{Paths in trees}

The analogous question of Csikv\'ari's theorem can be asked for paths rather than walks. Trivially, in a tree $T$ on $n$ vertices every pair of vertices determines exactly one path. So there are $\binom{n}{2}$ non-trivial paths in total, in particular there are at most that many paths of a given length $\ell$. We would like to determine an exact upper bound and classify the extremal cases.

It turns out that the answer depends heavily on the parity of $\ell$. If $\ell$ is odd, then the graph induced on $V(T)$ by the paths of length $\ell$ is bipartite, since $T$ is bipartite and paths of odd length connect vertices in different parts. So the obvious estimate gives an upper bound of $n^2/4$, which with a little thought can be improved to $n(n-\ell+1)/4$. On the other hand the construction of a double broom, i.e.~a path of length $\ell-2$ with equally many vertices attached to its ends yields $\lfloor(n-\ell+1)/2\rfloor\cdot \lceil(n-\ell+1)/2\rceil\ $ paths of length $\ell$. We shall prove that this construction is indeed optimal.

For even values of $\ell$ there is a better construction, namely a $p$-broom, i.e. a central vertex $v_0$ with $p$ many brooms attached to to it, where each broom consists of a path of length $(\ell-2)/2$ with several edges attached to the opposite end of $v_0$. Note that for $p=2$ this definition is consistent with the above notion of a double broom, whereas for $p>2$ it makes only sense if $\ell$ is even. For $\ell=2$ there is only one $p$-broom, namely the $n-1$-broom, or simply a star. Clearly, the star provides the optimal example for $\ell=2$, as every two edges in a star are adjacent and thus define a $2$-path. As in Tur\'an's theorem, for a given $p$ the maximal number of paths in a $p$-broom is attained if the sizes of the brooms are as equal as possible.

We prove that the maximal number of $\ell$-paths is always realized by some $p$-broom.

Let us say that vertices $v_1$ and $v_2$ of $T$ are $\ell$-neighbours if $d_T(v_1,v_2)=\ell$, where $d_T$ denotes the usual distance function on $T$. Let the \emph{valency} of $v$, in notation $r(v)$, be the number of its $\ell$-neighbours. Finally, let $R(T)$ be the total number of paths of length $\ell$ in $T$.

\begin{theorem}\label{thm: treebrooms}
For every $\ell$ there is a $p$ such that the maximal number of $\ell$-paths in a tree on $n$ vertices is attained for a $p$-broom.
\end{theorem}

\textbf{Proof.} Suppose that the tree $T$ realizes the maximum. Observe that if $T$ has two leaves $v$ and $w$ of distance other than $\ell$, then we must have $r(v)=r(w)$. Indeed, otherwise we could remove the leaf with the smaller valency, say $v$, and add a clone of $w$ instead. Let us call this operation the $DC$-transformation (as in delete-clone).

Since $r(v)=r(w)$ for any two leaves $v$ and $w$ of $d(v,w)\neq \ell$, we can freely apply the $DC$-transformation to any such pair $(v,w)$. Suppose now $d(v,w)>\ell$. Define $v'$ to be the unique vertex lying between $v$ and $w$ such that $d(v,v')=\ell$. Let $W$ be the set of all vertices that are separated from $v$ by $v'$; obviously, $w\in W$. By applying consecutive $DC$-transformations we can replace $W$ with clones of $v$. The resulting tree $T'$ would still achieve the maximal possible $R$, but the number of leaves would go up (since we create a clone-leaf for each vertex in $W$ and $v'$ is a new leaf).

So taking $T$ to be the tree with the maximal number of leaves amongst all $n$-vertex trees of maximal $R$, we can assume that $diam(T) = \ell$. Suppose now that there are two leaves $v$ and $w$ such that $2<d(v,w)<\ell$. Applying consecutive $DC$-transformation to $v$ and his clones on the one side and $w$ on the other side we can decrease the number of clone-classes of leaves, without decreasing $R$ or the total number of leaves. So we may assume that no such $v$ and $w$ exist.

To summarize the above arguments, we may assume that $T$ is of diameter $\ell$ and that any two leaves lie at distance either $2$ or $\ell$ from each other. Therefore $T$ comprises a number of stars, whose centres lie at mutual distance $\ell-2$. This is only possible if $T$ is a $p$-broom for some $p$; if $\ell$ is odd, $p$ must be equal $2$.

We actually have shown more, namely that every extremal tree $T$ can be obtained from a $p$-broom by applying a series of $R$-preserving inverse $DC$-transformations. It follows that the brooms are essentially unique extremal examples. The only exceptional cases occur when $n=\ell+2$, in which case the path of length $\ell+1$ is just as good as the broom, and when $\ell=3$, where there are more extremal examples, which all have a diameter of at most $4$ and are easy to classify. To be precise, every such example consists of a central vertex $v_0$ which has $a$ neighbours and $b$ vertices of distance $2$ from $v_0$. Such a tree has $b(a-1)$ paths of length three, which equals $\lfloor(n-2)^2/4\rfloor$ whenever $b = a-1 = (n-2)/2$.

Finally, let us return to the case when $\ell\ge 4$ is even and discuss briefly for what value of $p$ is the number of $\ell$-paths in a $p$-broom maximal. As was mentioned above, the numbers of leaves in the brooms must be equidistributed. As in Tur\'an's theorem, it follows that $f(p)$, the number of paths in such a $p$-broom satisfies \[\binom{p}{2}\left(\frac{n-1}{p}-\frac{\ell-2}{2}\right)^2-\frac{p}{8}\leq f(p)\leq \binom{p}{2}\left(\frac{n-1}{p}-\frac{\ell-2}{2}\right)^2.\] After some straightforward calculations involving differentiation and solving a quadratic equation, it turns out that the above error term of at most $p/8$ can be ignored and the maximum of the function $f$ lies within $1$ from \[p_{opt}= \frac{1}{4}+\sqrt{\frac{1}{16}+\frac{n-1}{\ell-2}}.\]

So, for example if we fix $\ell$ and let $n$ go to infinity, we obtain that $p$ is about $\sqrt{n/(\ell-2)}$ and \[f(p) = \binom{n}{2}-(\ell-2)^{1/2}n^{3/2}+O(n).\] 

If $p_{opt}$ is an integer, let's say $p_{opt}=k$, then $n-1=(\ell-2)(k^2-k/2)$, the maximum of $f$ is attained at $p_{opt}$ and equals \[f(p_{opt})=f(k)=\binom{k}{2}(\ell-2)^2(k-1)^2.\]

\end{document}